\documentclass[12pt]{article}
\usepackage[utf8]{inputenc}
\usepackage[a4paper, margin = 0.5in, bottom = 1in]{geometry}
\usepackage{soul}
\usepackage{authblk}

\date{}
\usepackage{amsmath}
\usepackage{amsfonts}
\usepackage{amssymb}
\usepackage{mathrsfs}
\usepackage{bm}
\usepackage{multicol}
\usepackage{wasysym}
\usepackage{graphicx}
\usepackage{caption}
\usepackage{cite}
\usepackage{subcaption}
\usepackage{skak}
\usepackage[table]{xcolor}
\usepackage{bbm}
\usepackage{xfrac}

\usepackage{amsfonts}
\usepackage{subcaption}
\usepackage{lmodern}
\usepackage{xcolor}
\usepackage{csquotes}
\usepackage{mathtools}
\usepackage{float}
\usepackage{amsmath,amssymb,amsthm,comment}
\usepackage{amsfonts}
\usepackage{subcaption}
\usepackage{xcolor}
\usepackage{csquotes}
\usepackage{mathtools}
\usepackage{float}
\usepackage{amsmath,amssymb,amsthm}
\usepackage[colorlinks,citecolor=red,urlcolor=blue,linkcolor=blue]{hyperref}
\usepackage{url}
\usepackage{mathrsfs,dsfont}
\usepackage{graphicx}
\usepackage{enumerate}
\usepackage{nicefrac}
\usepackage[font={normalsize},width=1.12\linewidth]{caption}
\usepackage{comment}
\usepackage{hyperref}
\usepackage{pgfplots}
\pgfplotsset{compat=1.15}
\usepackage{mathrsfs}
\usetikzlibrary{arrows}

\newtheorem{nummer}{ }
\newtheorem{thm}[nummer]{Theorem}

\newtheorem{lem}[nummer]{Lemma}
\newtheorem{cor}[nummer]{Corollary}

\newtheorem{defi}[nummer]{Definition}

\newcommand{\Mod}[1]{(\mathrm{mod}\; #1)}

\newtheorem{citedthm}{Theorem}

\makeatletter
\def\opargproof[#1]{\par\noindent {\bf #1 }}
 \makeatother

\begin{document}
\medskip\medskip
\begin{center}
\vspace*{50pt}
{\LARGE\bf The $m$-th Element of a Sidon Set}

\bigskip
{\small 
R. Balasubramanian}\\[1.2ex] 
{\scriptsize  The Institute of Mathematical Sciences, Chennai  \\ balu@imsc.res.in}\\[1.8ex]

{\small Sayan Dutta}\\[1.2ex] 
{\scriptsize Département de mathématiques et de statistique, Université de Montréal\\
sayan.dutta@umontreal.ca\\
\href{https://sites.google.com/view/sayan-dutta-homepage}{https://sites.google.com/view/sayan-dutta-homepage}}
\\[1.8ex]

\end{center}

\begin{abstract}
We prove that if $A=\{a_1,\dots ,a_{|A|}\}\subset \{1,2,\dots ,n\}$ is a Sidon set so that $|A|=n^{1/2}-L^\prime$, then
$$a_m = m\cdot n^{1/2} + \mathcal O\left( n^{7/8}\right) + \mathcal O\left(L^{1/2}\cdot n^{3/4}\right)$$
where $L=\max\{0,L^\prime\}$. As an application of this, we give easy proofs of some previously derived results. We proceed on to proving that for a dense Sidon set $S$ and for any $\varepsilon >0$, we have
$$\sum_{a\in S} a = \frac 12 n^{3/2} + \mathcal O \left (n^{11/8} \right )$$
for all $n\le N$ but at most $\mathcal O_{\varepsilon} \left (N^{\frac 45 + \varepsilon} \right )$ exceptions. 
\end{abstract}

\section{Introduction}
A set of positive integers $A\subset \mathbb N$ is called a \textit{Sidon Set} or a \textit{Sidon Sequence} if the equation $a+b=c+d$ does not have any non-trivial solutions in $A$. They were named after Hungarian mathematician Simon Sidon who was inspired by certain problems in Fourier series to ask Erd\H{o}s about the possible growth of such sequences.

Since then, there has been an extensive amount of literature on this topic exploring a plethora of different questions about finite and infinite Sidon sets \cite{lit}. Addressing the original question of Sidon, Erd\H{o}s conjectured \cite{conj} (and offered \$\,500 for a proof or disproof) that if $S(n)$ denotes the maximum possible cardinality of a Sidon subset of $[n] := \{1,2,\dots ,n\}$, then
$$S(n) < n^{1/2} + o(n^\varepsilon)$$
for all $\varepsilon >0$.

Several attempts on this problem by several different authors have only yielded
$$S(n)\le n^{1/2} + \mathcal O \left(n^{1/4}\right)$$
although a variety of different techniques have been tried \cite{lind, ruz, cil, Gra}. A classic result by Bose and Chowla \cite{bc} shows that the conjectured bound cannot be improved.

The most recent result on this was provided by Balogh, Füredi and Roy \cite{bal} proving
\begin{equation}\label{bal}
    S(n)< n^{1/2} + (1-\gamma) n^{1/4}
\end{equation}
for some $\gamma\ge 0.002$. In the same spirit, O'Bryant \cite{kev} improved the constant from $0.998$ to $0.99703$, and then again to $0.98183$ with Carter and Hunter \cite{dzk}.

Another recent improvement in the theory of Sidon sets is a solution to the Erdős Sárközy Sós problem on asymptotic Sidon bases of order $3$. After a series of developments made by Deshoulliers and Plagne \cite{dp}, Kiss \cite{kis}, Kiss, Rozgonyi and Sándor \cite{krs} and Cilleruelo \cite{cilced}, it was finally solved by Pilatte \cite{ced} using some recent results of Sawin \cite{saw}.

However our focus is on finding a formula for the $m$-th element of a finite Sidon set. As an application of our formula, we will indicate easier proofs of the theorems recently obtained by Ding \cite{ding1, ding2}. 

Our main tool is a theorem due to Cilleruelo.
\begin{citedthm}[Cilleruelo \cite{main}]\label{main}
   Let $A\subset [n]$ be a Sidon set with $\left\lvert A\right\rvert = n^{1/2}-L^\prime$. Then, every subinterval $I\subset [1,n]$ with length $cn$ contains $c\left\lvert A\right\rvert + E_I$ elements of $A$ where
   $$\left \lvert E_I \right \rvert \le 52\; n^{1/4} \left(1 + c^{1/2}\cdot n^{1/8}\right)\left(1+L^{1/2}\cdot n^{-1/8}\right)$$
   where $L = \max \{0,L^\prime\}$.
\end{citedthm}

\section{Dense Sidon Sets}
\begin{defi}
    A  Sidon subset $A\subset [n]$ is called \textit{dense} if
    $$\left\lvert A\right\rvert = \max \left\lvert S\right\rvert$$
    where the maximum is taken over all Sidon subsets of $[n]$.
\end{defi}

The structure of dense Sidon sets has a rich literature \cite{dense, sean} and classic constructions by Erd\H{o}s–Turán \cite{ertur}, Singer \cite{sin}, Bose \cite{bose}, Spence \cite{gan, ruz}, Hughes \cite{hug} and Cilleruelo \cite{cildense} have established that a dense Sidon set $A$ satisfies $\left\lvert A\right\rvert \ge \left(1-o(1)\right)\sqrt{n}$. As remarked by Ruzsa, ``somehow all known constructions of dense Sidon sets involve the primes" \cite{ruzdense}.

We will begin by proving a lemma that will justify an assumption we will make in the results to follow. This requires a classic result of Bose and a recent result of Baker, Harman and Pintz. 

\begin{citedthm}[R. C. Bose \cite{bose}]\label{bose}
    For a prime $p$, there are at least $p$ elements in $[p^2-1]$ such that the sums of two of these elements are different modulo $p^2-1$.
\end{citedthm}
\begin{citedthm}[R.C. Baker, G. Harman, J. Pintz. \cite{bhp}]\label{bhp}
    We have
    $$p_{k+1}-p_k \ll p_k^{0.525}$$
    where $p_k$ is the $k$-th prime.
\end{citedthm}

\begin{lem}\label{L_n}
    Let $A\subset [n]$ be a dense Sidon set so that $|A|=n^{1/2}-L$. Then $-n^{1/4}\le L\ll n^{21/80}$.
\end{lem}
\begin{proof}
    As discussed earlier, from Equation (\ref{bal}), we immediately have $L \ge -n^{1/4}$.
    
    For the next part, let us denote
    $$S_n = \max \left\lvert S\right\rvert$$
    over all Sidon subsets $S\subset [n]$.

    Also assume
    $$p_k^2 -1 < n \le p_{k+1}^2 -1$$
    where $p_i$ is the $i$-th prime. This gives
    $$p_k\le S_n$$
    using Theorem \ref{bose}.

    But from Theorem \ref{bhp}, we have
    $$p_{k+1}-p_k \ll p_k^{21/40} \le n^{21/80}$$
    hence giving
    $$p_k \ge \sqrt n - \mathcal O\left (n^{21/80}\right)$$
    so that
    $$\sqrt n - L = |A| = S_n \ge \sqrt n - \mathcal O\left (n^{21/80}\right)$$
    thus completing the proof.
\end{proof}

\textit{Remark}: Notice that if we assume $p_{k+1}-p_k \ll \sqrt{p_k}$, then an exact same line of argument will give $L\ll n^{1/4}$.

\section{The $m$-th element}
As advertised, we will now give the main theorem of the paper.
\begin{thm}\label{a_m} 
    Let $A = \{a_1,\dots , a_{\left\lvert A\right\rvert}\}\subset [n]$ be a dense Sidon set with $|A|=n^{1/2}-L^\prime$. Then
    $$a_m = m\cdot n^{1/2} + \mathcal O\left( n^{7/8}\right) + \mathcal O\left(L^{1/2}\cdot n^{3/4}\right)$$
    where $L=\max\{0,L^\prime\}$.
\end{thm}
\begin{proof}
    Consider a Sidon set $A=\{a_1,\dots ,a_{|A|}\}\subset [n]$ so that $|A|=n^{1/2}-L$ with $L\le n^{\frac {21}{80}}$. Let $I=[a_m]\subset [n]$ . Let $a_m=cn$ so that $c=\dfrac{a_m}{n}\le 1$. So, $|I|=cn$ and $\left\lvert A\cap [a_m]\right \rvert = m$.

    By Theorem \ref{main}, the number of elements in an interval $I$ of length $cn$ is
    $$c\left\lvert A\right\rvert + \mathcal O\left(n^{1/4}\right) + \mathcal O\left(n^{3/8}\right) + \mathcal O\left(L^{1/2} \cdot  n^{1/8}\right) + \mathcal O\left(L^{1/2} \cdot n^{1/4}\right)$$
    for $c\le 1$.

    This gives
    \begin{align*}
        m
        &= c\left( n^{1/2} +\mathcal O\left(L\right)\right) + \mathcal O\left( n^{3/8}\right) + \mathcal O\left(L^{1/2}\cdot n^{1/4}\right)\\
        &= c n^{1/2} +\mathcal O\left(L\right) + \mathcal O\left( n^{3/8}\right) + \mathcal O\left(L^{1/2}\cdot n^{1/4}\right)
    \end{align*}
    and hence putting $c=\dfrac{a_m}{n}$, multiplying by $n^{1/2}$, and rearranging, we get
    $$a_m = m\cdot n^{1/2} + \mathcal O\left( n^{7/8}\right) + \mathcal O\left(L^{1/2}\cdot n^{3/4}\right)$$
    thus completing the proof.
\end{proof}

This gives us the immediate corollary.
\begin{cor}\label{form}
    Let $A=\{a_1,\dots ,a_{|A|}\}\subset [n]$ be a Sidon set so that $|A|=n^{1/2}-L^\prime$ with $L\le n^{\frac {21}{80}}$. Then,
    $$a_m^\ell = m^\ell \cdot n^{\frac \ell 2} + \mathcal O \left( m^{\ell -1}\cdot n^{\frac{4\ell +3}{8}} \right) + \mathcal O\left( m^{\ell -1}\cdot L^{1/2}\cdot n^{\frac{2\ell +1}{4}}\right) +\mathcal O \left(n^{\frac{7\ell}8}\right) + \mathcal O\left(L^{\frac{\ell}{2}}\cdot n^{\frac{3\ell}{4}}\right)$$
    for any positive integer $\ell$.
\end{cor}
\begin{proof}
    This follows from the more general statement that if
    $$f(x) = m(x) + \mathcal O \left(e(x)\right)$$
    then
    $$\left(f(x)\right)^\ell = \left(m(x)\right)^\ell+\mathcal O\left(e(x)\left(m(x)\right)^{\ell-1}\right) + \mathcal O \left(e(x)^\ell\right)$$
    using the Binomial Theorem.
\end{proof}

We now also have the following corollaries that were also derived in \cite{ding1} and \cite{ding2}.

\begin{cor}\label{suma}
    Let $A=\{a_1,\dots ,a_{|A|}\}\subset [n]$ be a Sidon set so that $|A|=n^{1/2}-L^\prime$. Then,
    $$\sum_{a\in A} a = \frac 12 \cdot n^{3/2} + \mathcal O \left(n^{\frac{11}{8}}\right) + \mathcal O \left( L^{1/2}\cdot n^{5/4}\right)$$
    for $L\le n^{\frac {21}{80}}$.
\end{cor}
\begin{proof}
    Using Theorem \ref{a_m} and using the fact that $\left\lvert A\right\rvert < 2\sqrt n$, we have
    \begin{align*}
        \sum_{m=1}^{\left\lvert A\right\rvert} a_m &= \frac{\left\lvert A\right\rvert \left(\left\lvert A\right\rvert +1\right)}{2}\sqrt{n} + \mathcal{O} \left(n^{7/8} \left\lvert A\right\rvert\right) + \mathcal{O}\left(L^{1/2} n^{3/4} \left\lvert A\right\rvert\right)\\
        &= \frac{\sqrt n}{2}\left(\sqrt n - L\right)^2 + \frac{\sqrt n}{2}\left(\sqrt n - L\right) + \mathcal{O}\left(n^{7/8}\sqrt n\right) + \mathcal{O} \left(L^{1/2}\cdot n^{3/4} \sqrt n\right)\\
        &= \frac{\sqrt n}{2}\left(n + L^2 -2L\sqrt n\right) + \frac{\sqrt n}{2}\left(\sqrt n - L\right) + \mathcal O \left(n^{\frac{11}{8}}\right) + \mathcal O\left(L n^{7/8}\right) + \mathcal O \left(L^{1/2} n^{5/4}\right) + \mathcal O\left(L^{3/2} n^{3/4}\right)\\
        &= \frac 12 \cdot n^{3/2} + \mathcal O \left(n^{\frac{11}{8}}\right) + \mathcal O \left( L^{1/2}\cdot n^{5/4}\right)
    \end{align*}
    hence completing the proof.
\end{proof}

\begin{cor}\label{sumal}
    Let $A=\{a_1,\dots ,a_{|A|}\}\subset [n]$ be a Sidon set so that $|A|=n^{1/2}-L^\prime$. Then,
    $$\sum_{a\in A} a^\ell = \frac {1}{\ell +1} \cdot n^{\frac{2\ell +1}{2}} + \mathcal O \left( n^{\frac{8\ell +3}{8}} \right) + \mathcal O\left( L^{1/2}\cdot n^{\frac{4\ell +1}{4}}\right)$$
    for $L\le n^{\frac {21}{80}}$.
\end{cor}
\begin{proof}
    Using Corollary \ref{form} in exactly the same way as in the previous proof, we get
    \begin{align*}
        \sum_{m=1}^{\left\lvert A\right\rvert} a_m^\ell &= n^{\frac{\ell}{2}}\cdot \sum_{m=1}^{\left\lvert A\right\rvert} m^\ell + \mathcal O\left(n^{\frac{4\ell+3}{8}}\cdot \sum_{m=1}^{\left\lvert A\right\rvert} m^{\ell-1}\right) + \mathcal O\left(L^{1/2}\cdot n^{\frac{2\ell +1}{4}}\cdot\sum_{m=1}^{\left\lvert A\right\rvert} m^{\ell -1}\right)\\
        &\quad \quad \quad \quad \quad \quad \, + \mathcal O \left(n^{\frac{7\ell}8}\cdot \sqrt n\right) + \mathcal O\left(L^{\frac{\ell}{2}}\cdot n^{\frac{3\ell}{4}}\cdot \sqrt n\right)\\
        &= \frac{1}{\ell +1} \left(\sqrt n\right)^{\ell +1}\cdot n^{\frac{\ell}{2}} + \mathcal O \left(\left(\sqrt n\right)^\ell\cdot n^{\frac{4\ell +3}{8}}\right) + \mathcal O \left(\left(\sqrt n\right)^{\ell} \cdot L^{1/2} \cdot n^{\frac{2\ell +1}{4}}\right)\\
        &= \frac {1}{\ell +1} \cdot n^{\frac{2\ell +1}{2}} + \mathcal O \left( n^{\frac{8\ell +3}{8}} \right) + \mathcal O\left( L^{1/2}\cdot n^{\frac{4\ell +1}{4}}\right)
    \end{align*}
    hence completing the proof.
\end{proof}

\textit{Remark}: It should be noted that Theorem \ref{a_m} and Corollary \ref{form} are valid formulas for all $m$. However, they are only useful as an asymptotic formula when $m$ is close to $\sqrt n$.

\section{Improvements}
We will now provide an improvement to Theorem 1.4 and Corollary 1.5 of Ding \cite{ding2}. We will require the following result due to Heath-Brown which is an improvement on a series of similar results by Wolke \cite{wol}, Heath-Brown \cite{hb1, hb2}, Peck \cite{peck} and Matomäki \cite{ma}.
\begin{citedthm}[R. Heath-Brown \cite{hb}]\label{hb}
    For any $\varepsilon>0$ we have
    $$\sum_{\substack{p_n\le x\\ p_{n+1}-p_n\ge \sqrt{p_n}}} \left(p_{n+1}-p_n\right) \ll_{\varepsilon} x^{\frac 35 + \varepsilon}$$
    where $p_n$ is the $n$-th prime.
\end{citedthm}

We finally express a more quantitative version of the fact that the third error term can almost always be dropped.
\begin{thm}\label{improv}
    Let $S$ be a dense Sidon set in $[n]$. Then, for any $\varepsilon >0$, we have
    $$\sum_{a\in S} a = \frac 12 n^{3/2} + \mathcal O \left (n^{\frac {11}{8}} \right )$$
    for all $n\le N$ but at most $\mathcal O_{\varepsilon} \left (N^{\frac 45 + \varepsilon} \right )$ exceptions.
\end{thm}
\begin{proof}
We have already proven that
$$\sum_{a\in S} a = \frac 12 \cdot n^{3/2} + \mathcal O \left(n^{\frac{11}{8}}\right) + \mathcal O \left( L^{1/2}\cdot n^{5/4}\right)$$
for any choice of $L$.

So, it is enough to show that $\left\lvert L\right\rvert \ll n^{1/4}$ happens almost always with a small number of exceptions. As noticed in the remark under Lemma \ref{L_n}, denoting $p_m$ to be the $m$-th prime, if we have
$$p_{m+1}-p_m < \sqrt{p_m}$$
then for all $n\in \left [p_m^2-1, p_{m+1}^2-1 \right]$, we have $\left\lvert L\right\rvert\ll n^{1/4}$.

So, the set of exceptions are among those $n\in \left [ p_m^2-1, p_{m+1}^2-1\right ]$ with $p_{m+1}-p_m\ge \sqrt{p_m}$. So, the number of exceptions is bounded by
\begin{align*}
    \sum_{\substack{p_n\le \sqrt{N+1}\\ p_{n+1}-p_n \ge \sqrt{p_n}}} \left(p_{n+1}^2 - p_n^2\right)
    &= \sum_{\substack{p_n\le \sqrt{N+1}\\ p_{n+1}-p_n \ge \sqrt{p_n}}} \left(p_{n+1} + p_n\right)\left(p_{n+1} - p_n\right)\\\\
    &\ll \sqrt N \cdot \sum_{\substack{p_n\le \sqrt{N+1}\\ p_{n+1}-p_n \ge \sqrt{p_n}}} \left(p_{n+1} - p_n\right)\\\\
    &\ll_{\varepsilon} \sqrt N \cdot N^{\frac 3{10} + \varepsilon} = N^{\frac 45 + \varepsilon}
\end{align*}
where the first inequality uses the fact that $p_{n+1}+p_n\le 3\sqrt{N+1}$ and the second inequality follows from Theorem \ref{hb}.

This completes our proof.
\end{proof}

\textit{Remark}: An exact same argument gives
$$\sum_{a\in A} a^\ell = \frac {1}{\ell +1} \cdot n^{\frac{2\ell +1}{2}} + \mathcal O \left( n^{\frac{8\ell +3}{8}} \right)$$
for all $n\le N$ but at most $\mathcal O_{\varepsilon} \left (N^{\frac 45 + \varepsilon} \right )$ exceptions.

We conclude by making an observation on Theorem 1.2 of Ding \cite{ding1}. It essentially says that the sum of elements of a Sidon set $A$ in the residue class $i\; \Mod m$ is asymptotically $\frac{1}{2m} n^{3/2}$ provided that there is a function $f(n)\to \infty$ as $n\to \infty$ such that for all $t\in \left(\frac{n}{f(n)},n\right)$, we have $A(t):=A\cap (0,t)>\sqrt t$. However, it seems that this condition will rarely be satisfied for a dense Sidon set. For example, let us take $t=0.01 n$. By Theorem \ref{main}, $A(t)\sim 0.01 \sqrt n < 0.1 \sqrt n = \sqrt t$. In other words, Theorem \ref{main} prevents the said condition from being true under the assumption that $A(n)>\sqrt n$.

It should also be noted that the recent result given in \cite{li} improves Lemma \ref{L_n} to $L\ll n^{0.26}$. Similarly, the improvement given in \cite{olli} proves that the number of exceptions in Theorem \ref{improv} is $\mathcal O_{\varepsilon} \left (N^{0.785+\varepsilon} \right )$.

\section{Acknowledgements}
We would like to thank Prof. Kevin O'Bryant, Prof. Kaisa Matomäki and Runbo Li for bringing \cite{kev}, \cite{dzk}, \cite{olli} and \cite{li} to our attention. We would also like to thank the anonymous referee for his/her useful comments.

\bibliographystyle{plain}
{\footnotesize
\bibliography{Sidon}}

\begin{thebibliography}{10}

\bibitem{bhp}
R.C. Baker, Glyn Harman, and J.~Pintz.
\newblock On the difference between consecutive primes {II}.
\newblock {\em Proceedings of the London Mathematical Society}, 83:532--562, 2001.
\newblock \url{https://londmathsoc.onlinelibrary.wiley.com/doi/abs/10.1112/plms/83.3.532}.

\bibitem{bal}
József Balogh, Zoltán Füredi, and Souktik Roy.
\newblock An upper bound on the size of sidon sets.
\newblock {\em The American Mathematical Monthly}, 130(5):437--445, 2023.
\newblock \url{https://doi.org/10.1080/00029890.2023.2176667}.

\bibitem{bose}
R.~C. Bose.
\newblock An {A}ffine {A}nalogue of {S}inger's {T}heorem.
\newblock {\em {The Journal of the Indian Mathematical Society}}, 6:1--15, 1942.

\bibitem{bc}
R.C. Bose and S.~Chowla.
\newblock Theorems in the {A}dditive {T}heory of {N}umbers.
\newblock {\em Commentarii mathematici Helvetici}, 37:141--147, 1962/63.
\newblock \url{https://doi.org/10.1007/BF02566968}.

\bibitem{dzk}
Daniel Carter, Zach Hunter, and Kevin O'Bryant.
\newblock On the diameter of finite sidon sets.
\newblock {\em {On the diameter of finite Sidon sets}}, 175:108–126, 2025.
\newblock \url{https://link.springer.com/article/10.1007/s10474-024-01499-8}.

\bibitem{main}
Javier Cilleruelo.
\newblock Gaps in dense {S}idon sets.
\newblock {\em Integers}, 11(A11), 2000.
\newblock \url{https://matematicas.uam.es/~franciscojavier.cilleruelo/Papers/gap.pdf}.

\bibitem{cil}
Javier Cilleruelo.
\newblock Sidon sets in $\mathbb{N}^d$.
\newblock {\em Journal of Combinatorial Theory}, 117(7):857--871, 2010.
\newblock \url{https://doi.org/10.1016/j.jcta.2009.12.003}.

\bibitem{cildense}
Javier Cilleruelo.
\newblock Combinatorial problems in finite fields and {S}idon sets.
\newblock {\em Combinatorica}, 32:497--511, 2012.
\newblock \url{https://doi.org/10.1007/s00493-012-2819-4}.

\bibitem{cilced}
Javier Cilleruelo.
\newblock {On {S}idon sets and asymptotic bases}.
\newblock {\em Proceedings of the London Mathematical Society}, 111(5):1206--1230, 11 2015.
\newblock \url{https://academic.oup.com/plms/article-pdf/111/5/1206/6861604/pdv050.pdf}.

\bibitem{dp}
Jean-Marc Deshouillers and Alain Plagne.
\newblock A {S}idon basis.
\newblock {\em Acta Mathematica Hungarica}, 123(3):233–238, 2009.
\newblock \url{https://doi.org/10.1007/s10474-008-8097-3}.

\bibitem{ding1}
Yuchen Ding.
\newblock Sum of elements in finite {S}idon sets.
\newblock {\em International Journal of Number Theory}, 17(4):991--1001, 2021.
\newblock \url{https://doi.org/10.1142/S1793042121500196}.

\bibitem{ding2}
Yuchen Ding.
\newblock Sum of elements in finite {S}idon sets {II}.
\newblock {\em Publicationes Mathematicae Debrecen}, 103(1--2):243--256, 2023.
\newblock \url{https://doi.org/10.5486/PMD.2023.9595}.

\bibitem{dense}
Sean Eberhard and Freddie Manners.
\newblock The {A}pparent {S}tructure of {D}ense {S}idon {S}ets.
\newblock {\em The Electronic Journal of Combinatorics}, 30(P1.33), 2023.
\newblock \url{https://doi.org/10.37236/11191}.

\bibitem{conj}
Paul Erd\H{o}s.
\newblock Some problems in number theory, combinatorics and combinatorial geometry.
\newblock {\em Mathematica Pannonica}, 5:261–269, 1994.
\newblock \url{https://mathematica-pannonica.ttk.pte.hu/articles/mp05-2/mp05-2-261-269.pdf}.

\bibitem{ertur}
Paul Erdös and Pal Turán.
\newblock On a {P}roblem of {S}idon in {A}dditive {N}umber {T}heory, and on some {R}elated {P}roblems.
\newblock {\em Journal of the London Mathematical Society}, s1-16(4):212--215, 1941.
\newblock \url{https://doi.org/10.1112/jlms/s1-16.4.212}.

\bibitem{gan}
Michael~J. Ganley.
\newblock Direct product difference sets.
\newblock {\em Journal of Combinatorial Theory, Series A}, 23(3):321--332, 1977.
\newblock \url{https://doi.org/10.1016/0097-3165(77)90023-1}.

\bibitem{Gra}
S.~W. Graham.
\newblock {\em {$B_h$ sequences}}, pages {431--449}.
\newblock {Birkh{\"a}user Boston}, {Boston, MA}, 1996.
\newblock \url{10.1007/978-1-4612-4086-0_23}.

\bibitem{hb1}
Roger Heath-Brown.
\newblock The {D}ifferences between {C}onsecutive {P}rimes.
\newblock {\em Journal of the London Mathematical Society}, s2-18(1):7--13, 1978.
\newblock \url{https://doi.org/10.1112/jlms/s2-18.1.7}.

\bibitem{hb2}
Roger Heath-Brown.
\newblock The {D}ifferences between {C}onsecutive {P}rimes, {III}.
\newblock {\em Journal of the London Mathematical Society}, s2-20(2):177--178, 10 1979.
\newblock \url{10.1112/jlms/s2-20.2.177}.

\bibitem{hb}
Roger Heath-Brown.
\newblock {The {D}ifferences {B}etween {C}onsecutive {P}rimes, V}.
\newblock {\em International Mathematics Research Notices}, 2021(22):17514--17562, 12 2019.
\newblock \url{10.1093/imrn/rnz295}.

\bibitem{hug}
D.~R. Hughes.
\newblock Planar division neo-rings.
\newblock {\em Transactions of the American Mathematical Society}, 80:502--527, 1955.
\newblock \url{https://doi.org/10.1090/S0002-9947-1955-0073566-9}.

\bibitem{olli}
Olli Järviniemi.
\newblock On large differences between consecutive primes, 2022.
\newblock \url{https://arxiv.org/abs/2212.10965}.

\bibitem{kis}
S{\'a}ndor~Z. Kiss.
\newblock On {S}idon sets which are asymptotic bases.
\newblock {\em Acta Mathematica Hungarica}, 128:46 -- 58, 2010.
\newblock \url{https://doi.org/10.1007/s10474-010-9155-1}.

\bibitem{krs}
S{\'a}ndor~Z. Kiss, Eszter Rozgonyi, and Csaba S{\'a}ndor.
\newblock On {S}idon sets which are asymptotic bases of order $4$.
\newblock {\em Functiones et Approximatio Commentarii Mathematici}, 51(2):393 -- 413, 2014.
\newblock \url{https://doi.org/10.7169/facm/2014.51.2.10}.

\bibitem{li}
Runbo Li.
\newblock {The number of primes in short intervals and numerical calculations for Harman's sieve}, 2025.
\newblock \url{https://arxiv.org/abs/2308.04458}.

\bibitem{lind}
Bernt Lindström.
\newblock An inequality for ${B}_2$-sequences.
\newblock {\em Journal of Combinatorial Theory}, 6:211--212, 1969.
\newblock \url{https://doi.org/10.1016/S0021-9800(69)80124-9}.

\bibitem{ma}
Kaisa Matomaki.
\newblock Large {D}ifferences {B}etween {C}onsecutive {P}rimes.
\newblock {\em Quarterly Journal of Mathematics}, 58(4):489--518, aug 2007.
\newblock \url{10.1093/qmath/ham021}.

\bibitem{lit}
Kevin O'Bryant.
\newblock A complete annotated bibliography of work related to {S}idon sequences.
\newblock {\em Electronic Journal of Combinatorics}, DS11(39), 2004.
\newblock \url{https://doi.org/10.37236/32}.

\bibitem{kev}
Kevin O’Bryant.
\newblock On the size of finite sidon sets.
\newblock {\em Ukrains’kyi Matematychnyi Zhurnal}, 76(8):1192 -- 1206, September 2024.
\newblock \url{https://umj.imath.kiev.ua/index.php/umj/article/view/7858}.

\bibitem{peck}
A.~S. Peck.
\newblock Differences {B}etween {C}onsecutive {P}rimes.
\newblock {\em Proceedings of the London Mathematical Society}, 76(1):33--69, 1998.
\newblock \url{https://doi.org/10.1112/S0024611598000021}.

\bibitem{ced}
Cédric Pilatte.
\newblock A solution to the {E}rdős–{S}árközy–{S}ós problem on asymptotic {S}idon bases of order 3.
\newblock {\em Compositio Mathematica}, 160(6):1418--1432, 2024.
\newblock \url{https://doi:10.1112/S0010437X24007140}.

\bibitem{sean}
Sean Prendiville.
\newblock Solving equations in dense sidon sets.
\newblock {\em Mathematical Proceedings of the Cambridge Philosophical Society}, 173(1):25–34, 2022.
\newblock \url{10.1017/S0305004121000402}.

\bibitem{ruz}
Imre~Z. Ruzsa.
\newblock {Solving a linear equation in a set of integers I}.
\newblock {\em Acta Arithmetica}, 65(3):259--282, 1993.
\newblock \url{http://eudml.org/doc/206579}.

\bibitem{ruzdense}
Imre~Z. Ruzsa.
\newblock Erdős and the {I}ntegers.
\newblock {\em Journal of Number Theory}, 79(1):115--163, 1999.
\newblock \url{https://doi.org/10.1006/jnth.1999.2395}.

\bibitem{saw}
Will Sawin.
\newblock Square-root cancellation for sums of factorization functions over short intervals in function fields.
\newblock {\em Duke Mathematical Journal}, 170(5):997--1026, 2021.
\newblock \url{https://doi.org/10.1215/00127094-2020-0060}.

\bibitem{sin}
James~W. Singer.
\newblock A theorem in finite projective geometry and some applications to number theory.
\newblock {\em Transactions of the American Mathematical Society}, 43:377--385, 1938.
\newblock \url{https://api.semanticscholar.org/CorpusID:121112335}.

\bibitem{wol}
D.~Wolke.
\newblock Grosse {D}ifferenzen aufeinanderfolgender {P}rimzahlen.
\newblock {\em Mathematische Annalen}, 218:269--271, 1975.
\newblock \url{https://doi.org/10.1007/BF01349699}.

\end{thebibliography}

\end{document}